\documentclass[11pt]{amsart}
\usepackage{amsmath,amsfonts,amscd,amssymb,epsfig,euscript}
\usepackage{amsxtra,mathrsfs,bm}
\newtheorem{theorem}{Theorem}

\newtheorem{lemma}{Lemma}
\newtheorem{varthm}{Theorem}
\newtheorem*{Thm}{Theorem}
\newtheorem{corollary}{Corollary}
\theoremstyle{definition}

\theoremstyle{remark} \newtheorem{remark}[theorem]{Remark}
\numberwithin{equation}{section}
\newcommand{\field}[1]{\ensuremath{\mathbb{#1}}}
\newcommand{\CC}{\field{C}}

\newcommand{\HH}{\field{H}}

\newcommand{\RR}{\field{R}}

\newcommand{\ZZ}{\field{Z}}

\DeclareMathOperator{\im}{Im}
 
 \DeclareMathOperator{\re}{Re}
\DeclareMathOperator{\Tr}{Tr}
\newcommand{\del}{\partial}
\newcommand{\delb}{\bar\partial}

\newcommand{\C}{{\mathbb{C}}}

\newcommand{\curly}[1]{\mathscr{#1}}

\newcommand{\cE}{\curly{E}}
\newcommand{\cF}{\curly{F}}

\newcommand{\cH}{\curly{H}}

\newcommand{\s}{\gamma}

\newcommand{\bk}{\backslash}

\newcommand{\Gr}{\mathrm{SL}(2,\mathbb{R})}

\newcommand{\pa}{\partial}
\newcommand{\la}{\langle}
\newcommand{\ra}{\rangle}

\newcommand{\vep}{\varepsilon}

\newcommand{\z}{\bar{z}}

\DeclareMathOperator{\tr}{tr} \DeclareMathOperator{\Hom}{Hom}

\DeclareMathOperator{\ad}{ad}
 
 \DeclareMathOperator{\End}{End}
 \DeclareMathOperator{\Ad}{Ad}

\begin{document}
\title[Moduli of parabolic bundles]
{The first Chern form\\ on moduli of parabolic bundles}  
\date{\today}
\author{Leon A. Takhtajan}
\address{Department of Mathematics \\
Stony Brook University\\ Stony Brook, NY 11794-3651 \\ USA}
\email{leontak@math.sunysb.edu}
\author{Peter G. Zograf}
\address{Steklov Mathematical Institute \\
St. Petersburg, 191023, Russia}
\email{zograf@pdmi.ras.ru}

\begin{abstract}
For moduli space of stable parabolic bundles on a compact Riemann surface, we derive an explicit formula for 
the curvature of its canonical line bundle with respect to Quillen's metric and interprete it as a local index 
theorem for the family of $\delb$-operators in the associated parabolic endomorphism bundles. The formula consists 
of two terms: one standard (proportional to the canonical K\"ahler form on the moduli space), and one nonstandard, 
called a cuspidal defect, that is defined by means of special values of the Eisenstein-Maass series. The cuspidal 
defect is explicitly expressed through the curvature forms of certain natural line bundles on the moduli space 
related to the parabolic structure. We also compare our result with Witten's volume computation.
\end{abstract}

\maketitle

\tableofcontents

\section{Introduction}
Local index theorems for families of $\bar{\partial}$-operators provide local (i.e., valid on the level of 
differential forms) expressions for the Chern classes (forms) of the corresponding index bundles.
Historically first examples of such results belong to Quillen \cite{Q} and Belavin-Knizhnik \cite{BK}.
Quillen considered the case of Cauchy-Riemann, or $\bar{\partial}$-operators, in a vector bundle on
a Riemann surface. He observed that when the natural $L^2$-metric in the determinant index bundle
is divided by the determinant of the Laplace operator $\det\bar\partial^*\bar{\partial}$, its curvature
becomes proportional to a natural K\"ahler form on the parameter space. Belavin and Knizhnik
extended Quillen's result to the families of $\bar{\partial}$-operators on compact Riemann surfaces.
Both papers rely on heat kernel expansion techniques.

The pioneering work of Quillen and Belavin-Knizhnik initiated an extensive treatment of various forms
of local index theorems in the literature. E.g., in our papers \cite{TakZog87:Localindexthm} and 
\cite{Takhtajan-Zograf89} we rederived and refined these results using deformation theory (in particular, Teichm\"uller
theory). A similar approach was also used in \cite{Fay92}. Our technique proved to be applicable to the families of $\delb$-operators on punctured Riemann surfaces
\cite{Takhtajan-Zograf88, Takhtajan-Zograf91}, giving the first example of a local index theorem for 
families with non-compact fibres. In this case the spectrum of the Laplace operator contains an absolutely 
continuous part, so that the standard heat kernel definition of the regularized determinant (that enters Quillen's 
metric) is not applicable. Instead, we define it as a special value of the Selberg zeta function. The curvature 
(or the first Chern form) of the determinant bundle then splits into two terms: 
one being proportional to the Weil-Petersson K\"ahler form 
on the moduli space of punctured Riemann surfaces and the other being the K\"ahler form of a new K\"ahler
metric defined in terms of Eisenstein-Maass series (the so-called cuspidal defect arising from 
the absolutely continuous spectrum of the Laplacian). Details can be found in \cite{Takhtajan-Zograf88, Takhtajan-Zograf91}. Further refinements of these results were obtained by Weng et al. \cite{Weng,OTW} 
in terms of Deligne pairings and Arakelov geometry, and by Wolpert \cite{Wolpert} in terms
of complex differential geometry.

The present paper is a long overdue sequel to \cite{Takhtajan-Zograf89}. Here we combine the methods 
developed in \cite{Takhtajan-Zograf89} and \cite{Takhtajan-Zograf88,Takhtajan-Zograf91} to treat another 
example --- a local index theorem for the family of $\delb$-operators acting in the endomorphism bundles 
associated with the stable parabolic bundles on a compact Riemann surface.\footnote{For families of stable
parabolic bundles with nonzero rational weights a version of local index theorem was obtained in \cite{BR}.
In the situation considered in \cite{BR} the corresponding Laplace operators have no continuous spectrum,
and no cuspidal defect appears in this case; the resulting formula is very similar to the original Quillen's
one.} 

To be more precise, let $E$ be a stable parabolic vector bundle of rank $k$ on a compact Riemann surface $X$, with given weights and multiplicities at the marked points $P_{1},\dots, P_{n}$. According to the Mehta-Seshadri theorem \cite{Mehta-Seshadri80}, the bundle $E$ is associated with an irreducible unitary representation of the fundamental group of the non-compact Riemann surface $X_{0}=X\setminus\{P_{1},\dots,P_{n}\}$. Put $X_{0}\cong\Gamma\bk\HH$, where $\HH$ is the Poincar\'{e} model of the Lobatchevsky plane, and $\Gamma$ is a torsion-free Fuchsian group. 
Then there exists an irreducible representation
$\rho:\Gamma\rightarrow U(k)$ such that the spectrum of 
$\rho(S_{i})$, where $S_{i}$ is a parabolic generator of $\Gamma$ about the marked point $P_{i}$, $i=1,\dots,n$,  
is given by the exponents of weights at $P_{i}$, and such that $E\cong E^{\rho}$, where $E^{\rho}$ is a proper extension 
of the quotient bundle $E_{0}^{\rho}= \Gamma\bk(\HH\times\CC^{k})\rightarrow \Gamma\bk\HH\cong X_{0}$ to the compact surface
$X$. The Hermitian metric in the bundle $\End E_{0}^{\rho}$ (induced by the standard Hermitian metric in $\End\CC^{k}$) 
and the complete hyperbolic metric on $X_{0}$ define the Hodge $\ast$-operator in the vector spaces 
of $\End E_{0}$-valued $(p,q)$-forms on $X_{0}$. The Laplace operator in the bundle $\End E_{0}$ is defined by $\Delta=\delb^{\ast}\delb$, where 
$\delb=\delb_{\,\End E_{0}}$ and $\delb^{\ast}=-\ast\delb\ast$; it is a self-adjoint operator in the Hilbert space 
of $L^2$-sections of $\End E_{0}$ on $X_{0}$.
The isomorphism $\End E_{0}\cong \End E_{0}^{\rho}$ identifies
the Laplace operator $\Delta$ with the Laplace operator on $\HH$ acting in the space of $\End \CC^{k}$-valued functions on $\HH$
automorphic with respect to $\Gamma$ with the unitary representation $\Ad\rho$ that are square integrable on the
fundamental domain of $\Gamma$. We define its regularized determinant as 
$$\det\Delta=\left.\frac{\del}{\del s}\right|_{s=1}Z(s,\Gamma;\Ad\rho),$$
where $Z(s,\Gamma;\Ad\rho)$ is the Selberg  zeta function corresponding to $\Delta$ 
(see Section 2 for precise definitions and references).

The moduli space $\mathcal{N}$ of stable parabolic vector bundles of rank $k$ on $X$ is a complex manifold\footnote{More
precisely, a smooth quasiprojective variety. Its natural compactification --- the moduli space of {\em semistable} parabolic bundles --- is a normal projective variety that is {\em smooth} for a generic weight system; cf. \cite{Mehta-Seshadri80}.}.
The holomorphic tangent space $T_{\{E\}}\mathcal{N}$ 
at the point $\{E\}\in\mathcal{N}$ corresponding to the stable parabolic bundle $E$ is naturally isomorphic to the space $\cH^{0,1}(X_{0},\End E_{0})$ of square integrable harmonic $\End E_{0}$-valued $(0,1)$-forms on $X_{0}$.
The moduli space $\mathcal{N}$ carries a natural K\"{a}hler metric given by the Hodge inner product in the tangent spaces.
We denote by $\Vert\cdot\Vert^{2}$ the corresponding Hermitian metric in the canonical line bundle
$\lambda=\det\mathcal{N}$ (see Section \ref{moduli space}).

Our main result ---  Theorem  1 of Section 5 --- is an explicit computation of the curvature form of Quillen's metric
$\Vert\cdot\Vert_{Q}^{2}=\Vert\cdot\Vert^{2}(\det\Delta)^{-1}$ in the canonical line bundle $\lambda$ on $\mathcal{N}$. 
In addition to the term proportional to the K\"ahler form on $\mathcal{N}$, it contains an extra term, the so-called cuspidal
defect, that is due to the absolutely continuous spectrum of $\Delta$. It is explicitly defined in terms of the values at $s=1$
of the Eisenstein-Maass series for the group $\Gamma$
%that are associated with the absolutely continuous spectrum of 
%$\Delta$ 
%in the space of $\End\CC^{k}$-valued $(\Gamma,\, \Ad\rho)$-automorphic functions on $\HH$
(see Section 5 for the precise formulation). We also interpret the cuspidal defect in terms of the curvature forms of natural line
bundles on the moduli space $\mathcal{N}$ associated with the parabolic structures.
This result simplifies for the moduli space $\mathcal{N}_{0}$ of parabolic bundles with fixed determinant, well-defined when the parabolic structure is integral (see Corollary 1 of Section 5). In particular, it gives an alternative approach to computing volumes of moduli spaces of parabolic bundles. We also compare our computations with Witten's formula \cite{Witten}
for the symplectic volume of $\mathcal{N}$ in the simplest situation of a pointed torus. 

The content of the paper is as follows. In Sections 2 and 3 we collect the necessary facts about stable parabolic bundles on compact Riemann surfaces and their moduli spaces, as well as about the spectral theory of automorphic Laplacians. In Section 4 we derive the necessary variational formulas, introduce certain natural line bundles on moduli spaces and compute their curvatures. At last, in Section 5 we prove the main result --- Theorem 1.  

\noindent
\textbf{Acknowledgments.} The first author (LT) was partially supported by the NSF grants DMS-0204628
and DMS-0705263. The second author (PZ) was partially supported by the President of Russian Federation 
research grant NSh-U329.2006.1 and by the Russian Foundation for Basic Research grant 05-01-00899.

\section{Preliminaries}

\subsection{Parabolic bundles} \label{parabolic bundles}

Let $X$ be a compact Riemann surface of genus $g$ with a finite set $S=\{P_{1},\dots,P_{n}\}$ of marked points 
such that $2g+n-2>0$. According to \cite{Mehta-Seshadri80}, a holomorphic vector bundle $E$ on $X$ of rank $k$ is called a \emph{parabolic bundle} if it carries a \emph{parabolic structure} --- a flag $E_{P}=F_{1}E_{P}\supset F_{2}E_{P}\supset\dots\supset F_{r}E_{P}$ in the fibre $E_{P}$ and \emph{weights} $0\leq\alpha_{1}<\alpha_{2}<\dots<\alpha_{r}<1$ 
for each $P\in S$. The integers $k_{1}=\dim F_{1}E_{P}-\dim F_{2}E_{P},\dots,k_{r}=\dim F_{r}E_{P}$,
are called the \emph{multiplicities} of the weights $\alpha_{1},\dots,\alpha_{r}$. 
A morphism $f: E\rightarrow E'$ of parabolic vector bundles is a morphism of holomorphic vector bundles such that for every $P\in S$,
$f(F_{i}E_{P})\subset F_{j+1}E'_{P}$ whenever $\alpha_{i}>\alpha'_{j}$.
A subbundle $F\subset E$ is a \emph{parabolic subbundle} if for every $P\in S$ the parabolic structure in $F$ is a restriction of the parabolic structure in $E$.
The \emph{parabolic degree} of a parabolic bundle $E$ is defined as
$$\mathrm{par}\, \mathrm{deg}\, E=\mathrm{deg}\, E +\sum_{P\in S}\sum_{l=1}^{r(P)}k_{l}(P)\alpha_{l}(P),$$
where $\mathrm{deg}\, E$ is the degree of the underlying holomorphic vector bundle $E$. A parabolic bundle $E$ of parabolic
degree 0 is said to be stable \cite{Mehta-Seshadri80} if $\mathrm{par}\,\mathrm{deg}\,F<0$ for every parabolic subbundle $F$
of $E$. A theorem of Mehta-Seshadri \cite{Mehta-Seshadri80} generalizes the celebrated theorem of Narasimhan-Seshadri \cite{Narasimhan-Seshadri65} about stable vector bundles on a compact Riemann surface to the case of parabolic bundles. It states that stable parabolic bundles are precisely those associated with irreducible unitary representations of the fundamental group of the non-compact Riemann surface $X_{0}=X\setminus S$. 

A precise formulation is the following. 
By the uniformization theorem, $X_{0}\cong\Gamma\bk\HH$, where
$\HH=\{z=x+\sqrt{-1}\,y\in\CC\, |\, y>0\}$ is a Poincar\'{e} model of the Lobatchevsky (hyperbolic) plane, and $\Gamma$ is a torsion-free Fuchsian group generated by hyperbolic transformations $A_{1},B_{1},\dots,A_{g},B_{g}$ and parabolic transformations $S_{1},\dots,S_{n}$ satisfying the single relation 
$$A_{1}B_{1}A_{1}^{-1}B_{1}^{-1}\dots A_{g}B_{g}A_{g}^{-1}B_{g}^{-1}S_{1}\dots S_{n}=1.$$
Let $x_{1},\dots,x_{n}$ be the fixed points of the elements $S_{1},\dots,S_{n}$ (also called parabolic cusps),
and let $\overline\HH$ be the union of $\HH$ with the set of all parabolic cusps of $\Gamma$. 
There is a natural projection $\overline\HH\rightarrow\Gamma\bk\overline\HH$ such that $X\cong\Gamma\bk\overline\HH$.
The images of the cusps $x_{1},\dots,x_{n}\in\RR\cup\{\infty\}$ 
are the marked points $P_{1},\dots,P_{n}$ (see, e.g., \cite[Ch.~1]{Shimura94}). Let $\CC^{k}$ be the complex vector space
with the standard Hermitian inner product and the orthonormal basis, and let $U(k)$ be the group of $k\times k$ unitary
matrices. A unitary representation $\rho:\Gamma\rightarrow U(k)$ is called \emph{admissible} with respect to a given set 
of weights and multiplicities at $P_1,\dots,P_n$,   
if  for each $i=1,\dots, n$, we have $\rho(S_{i})=U_{i}D_{i}U_{i}^{-1}$
with unitary $U_{i}\in U(k)$ and diagonal $D_{i}=e^{2\pi\sqrt{-1}\mathcal{A}_{i}}$,
$\mathcal{A}_{i}=(\alpha^{i}_{1},\dots,\alpha^{i}_{r_{i}})$, where each $\alpha^{i}_{l}=\alpha_{l}(P_{i})$ is repeated $k^{i}_{l}=k_{l}(P_{i})$ times, $l=1,\dots,r_{i}=r(P_{i})$. Admissible matrices $\rho(S_{i})$ are parametrized by the flag varieties $\cF_{i}=U(k)/
U(k^{i}_{1})\times\dots\times U(k^{i}_{r_{i}})$, $i=1,\dots,n$. The group $\Gamma$ acts on the trivial bundle $\overline\HH\times\CC^{k}$ on $\overline\HH$ by the rule $(z,v)\mapsto (\gamma z,\rho(\gamma)v)$, where $z\in\overline\HH,\, v\in\CC^{k}$ and $\gamma\in\Gamma$. Take the sheaf of its bounded at the cusps $(\Gamma,\rho)$-invariant
sections. The direct image of this sheaf under the projection $\overline\HH\rightarrow X$ is a locally free sheaf of rank $k$
on $X$. The parabolic structure at the images of cusps is defined by the matrices $\rho(S_{i})$. This gives a parabolic vector bundle $E^{\rho}$ on the Riemann surface $X$ of parabolic degree $0$. Loosely speaking, the bundle $E^{\rho}$  
is an extension to $X$ of the quotient bundle $E_{0}^{\rho}\cong\Gamma\bk(\HH\times\CC^{k})\rightarrow\Gamma\bk\HH\cong X_{0}$. It is easy to describe the bundle $E^{\rho}$ explicitly in terms of the transition functions as in \cite[Remark 6.2]{Narasimhan-Seshadri65}.

\begin{Thm}[Mehta-Seshadri] A parabolic vector bundle 
$E$ of rank $k$ and $\mathrm{par}\, \mathrm{deg}\, E=0$ is stable if and only if it is isomorphic to a bundle $E^{\rho}$, where
$\rho: \Gamma\rightarrow U(k)$ is an irreducible representation of the group $\Gamma$ admissible  with respect to the set of weights and multiplicities of the parabolic structure of $E$. Moreover, stable parabolic bundles $E^{\rho_{1}}$ and $E^{\rho_{2}}$ are isomorphic if and only if representations $\rho_{1}$ and $\rho_{2}$ are equivalent.  
\end{Thm}

\begin{remark} The original proof in \cite{Mehta-Seshadri80} was of algebro-geometric nature and worked only for rational weight systems. Following Donaldson's ideas \cite{Donaldson83}, a more straightforward differential-geometric proof 
valid for arbitrary real weights was given in \cite{Biquard91}.
\end{remark}

The standard Hermitian metric in $\CC^{k}$ defines a $\Gamma$-invariant metric in the trivial
vector bundle $\overline\HH\times\CC^{k}\rightarrow \overline\HH$. It extends as a (pseudo)metric $h_{E}$ to the bundle $E=E^{\rho}$ that degenerates in the fibres over the points $P_{1},\dots,P_{n}$. Explicitly, choose $\sigma_{i}\in\Gr$ 
such that $\sigma_{i}\infty=x_{i}$ and $\sigma_{i}^{-1}S_{i}\sigma_{i}=\left(\begin{smallmatrix}1 & \pm 1\\ 0 & \;\;1\end{smallmatrix}\right)$,
and let $\zeta=e^{2\pi\sqrt{-1}\,\sigma^{-1}_{i}z}$ be a local coordinate
at $P_{i}\in X\cong\Gamma\bk\overline\HH$. 
 Then in terms of the 
trivialization of $E$ defined by $k$ local sections --- the columns of the matrix $U_{i}e^{2\pi\sqrt{-1}\zeta\mathcal{A}_{i}}$ --- the metric $h_{E}$ is given by the diagonal matrix $|\zeta|^{2\mathcal{A}_{i}}=(|\zeta|^{2\alpha_{1}^{i}},\dots,|\zeta|^{2\alpha^{i}_{r_{i}}})$, where each $|\zeta|^{2\alpha^{i}_{l}}$ is repeated $k^{i}_{l}$ times, $l=1,\dots, r_{i}$. The restriction of $h_{E}$ to the bundle $E_{0}=E_{0}^{\rho}$, which we denote by $h_{E_{0}}$, is non-degenerate.

\subsection{The endomorphism bundle} \label{EndE} 

Let $\mathrm{End}\,E_{0}$ be the bundle of endomorphisms  of the vector bundle $E_{0}$ on $X_{0}$. 
Its fibers have the structure of the Lie algebra $\mathfrak{gl}(k,\CC)$ with the bracket $[~,~]$ and the Killing form 
$\tr$. Together with the exterior multiplication in the space $C^{\bullet}(X_{0})$ of smooth differential forms on $X_{0}$, these operations induce the mappings
\begin{align*}
[~,~] &:C^{p}(X_{0},\mathrm{End}\,E_{0})\otimes C^{q}(X_{0},\mathrm{End}\,E_{0})\rightarrow C^{p+q}(X_{0},\mathrm{End}\,E_{0}) \\
\intertext{and}
\wedge & : C^{p}(X_{0},\mathrm{End}\,E_{0})\otimes C^{q}(X_{0},\mathrm{End}\,E_{0})\rightarrow C^{p+q}(X_{0}),
\end{align*}
where $C^{p}(X_{0},\mathrm{End}\,E_{0})$ is the space of smooth
$\mathrm{End}\,E_{0}$-valued $p$-forms on $X_{0}$. For $E_{0}=E_{0}^{\rho}$, the bundle  $\mathrm{End}\,E_{0}\cong\Gamma\bk(\HH\times\mathrm{End}\,\CC^{k})$ is the quotient bundle with respect to the
adjoint representation $\Ad\rho$ of the group $\Gamma$ in $\mathrm{End}\,\CC^{k}$. Explicitly, 
$\mathrm{Ad}\,\rho(\gamma) a=\rho(\gamma)a\rho(\gamma)^{-1}$, where $\gamma\in\Gamma$ and $a\in\mathrm{End}\,\CC^{k}$, and 
realized as a $k^{2}\times k^{2}$ matrix,
$\Ad\rho(\s)=(\rho(\s)\otimes(\rho(\s)^{-1})^{t})=(\rho(\s)\otimes\overline{\rho(\s)})$. In particular, if $u_{1},
u_{2}\in\CC^{k}$ are eigenvectors of $\rho(\s)$ with eigenvalues $e^{\sqrt{-1}\theta_{1}}$ and
$e^{\sqrt{-1}\theta_{2}}$, then $v=u_{1}\otimes \bar{u}_{2}\in\End\C^{k}$ --- a $k\times k$ matrix with the elements
$v_{lm}=u_{1l}\bar{u}_{2m}$ --- is an eigenvector for $\Ad\rho(\s)$ with eigenvalue $e^{\sqrt{-1}(\theta_{1}-\theta_{2})}$.

The Hermitian metric $h_{E_{0}}$ in the bundle $E_{0}$ naturally induces a Hermitian metric in the bundle $\End E_{0}$
that we denote by $h_{\mathrm{End}\,E_{0}}$.
The bundle $\End E_{0}$ has a canonical section --- the identity isomorphism $I$ of $E_{0}$, 
and decomposes into the orthogonal sum $\End E_{0}=\ad E_{0}\oplus\,\CC$
with respect to the metric $h_{\End E_{0}}$. Here $\ad E_{0}$ is the adjoint bundle (that is, the bundle of traceless
endomorphisms of $E_0$), and $\CC$ is understood as a trivial line bundle on $X_{0}$ spanned on the non-vanishing 
section $I$. 

\begin{remark}
Since $\End E_{0}$ is associated with the unitary representation $\Ad\rho$ of $\Gamma$, it can be extended to $X$
as a parabolic bundle $\End E = \End E^\rho$ of parabolic degree 0. However, we will not use this extension ---
we are going to work with $L^2$-sections of $\End E_{0}$ instead.
\end{remark}

\subsection{The Laplace operator}\label{Laplace}

Let $E$ be a stable parabolic vector bundle on $X$ and let $E_{0}$ be the restriction of $E$ to 
$X_{0}=X\setminus S\cong \Gamma\bk\HH$.
We use the hyperbolic (or Poincar\'{e}) metric on $X_{0}$ descended from $\HH$ and the Hermitian metric $h_{\End E_{0}}$ in $\End E_{0}$
to define the Hodge $\ast$-operator in the vector spaces $C^{p,q}(X_{0},\End E_{0})$ for $p,q=0,1$.
Let $C^{p,q}_{c}(X_{0},\End E_{0})$ denote the subspace of compactly supported $\mathrm{End}\,E_{0}$-valued $(p,q)$-forms on $X_{0}$. The completion of this space with respect to the Hodge inner product 
yields the Hilbert space $\mathfrak{H}^{p,q}(X_{0},\End E_{0})$. The Laplace operator in $C^{0,0}_{c}(X_{0},\End E_{0})$ is,
by definiton, $\Delta=\delb^{\ast}\delb$, where $\delb$ and its adjoint $\delb^{\ast}=-\ast\delb\ast$ 
are understood as operators from 
$C^{0,0}_{c}(X_{0},\End E_{0})$ to $C^{0,1}_{c}(X_{0},\End E_{0})$ and from $C^{0,1}_{c}(X_{0},\End E_{0})$ to $C^{0,0}_{c}(X_{0},\End E_{0})$ respectively. The Laplace operator admits a unique extension as a non-negative, self-adjoint operator in the Hilbert space $\mathfrak{H}^{0,0}(X_{0},\End E_{0})$, which we also denote by $\Delta$. Since the bundle $E$ 
is stable, the kernel $\ker\Delta=\ker\delb$ of the operator $\Delta$ is one-dimensional and is generated by the section $I$. Denote by $\mathfrak{H}_{0}^{0,0}(X_{0},\End E_{0})$ the orthogonal complement of $\ker\Delta$ in $\mathfrak{H}^{0,0}(X_{0},\End E_{0})$ and by
$\Delta_{0}$ --- the restriction of the operator $\Delta$ to  $\mathfrak{H}_{0}^{0,0}(X_{0},\End E_{0})$. Then $\ker\delb^{\ast}
=\cH^{0,1}(X_{0},\End E_{0})$ is the subspace of harmonic $(0,1)$-forms in the Hilbert space $\mathfrak{H}^{0,1}(X_{0},\End E_{0})$. The corresponding orthogonal projection
$$P: \mathfrak{H}^{0,1}(X_{0},\End E_{0})\rightarrow\cH^{0,1}(X_{0},\End E_{0})$$ is given by
$$P=I-\delb\Delta_{0}^{-1}\delb^{\ast},$$
where $I$ stands now for the identity operator in $\mathfrak{H}^{0,1}(X_{0},\End E_{0})$. 

When $E=E^{\rho}$, the spaces $C^{p,q}(X_{0},\End E_{0})$ are naturally identified with the vector spaces of smooth $\End\CC^{k}$-valued automorphic forms on $\HH$ with the transformation law
$$f(\gamma z)\gamma^{\prime}(z)^{p}\overline{\gamma^{\prime}(z)}^{q}=\Ad\rho(\gamma)f(z),\quad z\in\HH,\,\gamma\in\Gamma.$$
The Hodge inner product is then
\begin{equation} \label{inner-1}
\la f_{1},f_{2}\ra =\int\limits_{X}f_{1}\wedge\ast f_{2}=2^{p+q}\iint\limits_{F}\tr(f_{1}(z)f_{2}(z)^{\ast})y^{p+q-2}dxdy,
\end{equation}
where $f^{\ast}=\bar{f}^{\,t}$ is the Hermitian conjugate of $f\in\End\CC^{k}$, 
and $F$ is a fundamental domain for the group $\Gamma$ in $\HH$. The Hilbert spaces $\mathfrak{H}^{p,q}(X_{0},\End E_{0})$ are naturally identified with the Hilbert spaces  $\mathfrak{H}^{p,q}(\HH, \Gamma; \Ad\rho)$ of 
$(\Gamma, \Ad\rho)$-automorphic forms on $\HH$. 
We have 
$$\delb=\frac{\del}{\del\z}=\frac{1}{2}\left(\frac{\del}{\del x}+\sqrt{-1}\,\frac{\del}{\del y}\right),\quad\delb^{\ast}=-2y^{2}\frac{\del}{\del z}=-y^{2}\left(\frac{\del}{\del x}-\sqrt{-1}\,\frac{\del}{\del y}\right),$$ 
so that the Laplace operator has the form
$$\Delta=-2y^{2}\frac{\pa^{2}}{\pa z\pa\bar{z}}=-\frac{y^{2}}{2}\left(\frac{\del^{2}}{\del x^{2}}+\frac{\del^{2}}{\del y^{2}}\right).$$ 
The spectral decomposition of the Laplace operator in the Hilbert space $\mathfrak{H}^{0,0}(\HH, \Gamma; \Ad\rho)$
has been studied in detail in
\cite{Venkov81}\footnote{The operator usually considered in the spectral theory is $2\Delta$.}. The spectrum of $\Delta$ has both discrete and absolutely continuous parts. The latter covers the interval $\left[\frac{1}{8},\infty\right)$ 
with the multiplicity 
$$\sum_{i=1}^{n}\sum_{l=1}^{r_{i}}(k^{i}_{l})^{2}=nk^{2}-\sum_{i=1}^{n}\dim_{\RR} \cF_{i}.$$
The eigenfunctions of the continuous spectrum are given by the analytically continued Eisenstein-Maass series for the group
$\Gamma$ with the unitary representation $\Ad\rho$. To be more specific, consider the subspaces
$V_{i}=\ker (\Ad\rho(S_{i})-I)$ in $\End\CC^{k}$, $i=1,\dots,n$. The Eisenstein-Maass series corresponding to the cusp $x_{i}$ and a vector $v\in V_{i}$ is defined for $\re s>1$ by the following absolutely convergent series
$$E_{i}(z,v;s)=\sum_{\gamma\in\Gamma_{i}\bk\Gamma}\mathrm{Im}(\sigma_{i}^{-1}\gamma z)^{s}\Ad\rho(\gamma)^{-1}v,\quad i=1,\dots,n.$$
Here $\Gamma_{i}$ is the stabilizer of the cusp $x_{i}$ in $\Gamma$ --- the cyclic subgroup generated by $S_{i}$, and $\sigma_{i}\in\Gr$ is as in Section \ref{parabolic bundles}. Since $v^{\ast}\in V_{i}$ for $v\in V_{i}$ and 
the representation $\rho$ is unitary, $E_{i}(z,v;s)^{\ast}=E_{i}(z,v^{\ast},\bar{s})$.
The Eisenstein-Maass series $E_{i}(z,v;s)$  is $(\Gamma,\Ad\rho)$-automorphic --- that is, $E_{i}(\s z,v;s)=\Ad\rho(\s)E_{i}(z,v;s)$. It satisfies the differential equation
\begin{equation} \label{E-Laplace}
\Delta E_{i}(z,v;s)=\frac{s(1-s)}{2}E_{i}(z,v;s),
\end{equation}
and admits a meromorphic continuation to the whole complex $s$-plane. 
Since the representation $\rho: \Gamma\rightarrow U(k)$ is irreducible, it follows from \eqref{E-Laplace} that
for every $v\in V_{i}$ satisfying $\tr v=0$ the Eisenstein-Maass series $E_{i}(z,v;s)$ is regular at $s=1$.
Equation \eqref{E-Laplace} together with the property
$$E_{i}(\sigma_{j}(z+1),v;s)=\Ad\rho(S_{j})E_{i}(\sigma_{j}z,v;s)$$ 
yield the following asymptotic expansion of $E_{i}(z,v;1)$ with $\tr v=0$ at the cusps:
\begin{equation} \label{E-asymptotic} 
E_{i}(\sigma_{j}z,v;1)=\delta_{ij}\,y\cdot v + c_{ij}(v) + O(e^{-\pi y})\quad 
\text{as $y\rightarrow\infty$},
\end{equation}
where all $c_{ij}(v),\;i,j=1,\dots,n,$ belong to $\End\CC^{k}$ and satisfy $\tr c_{ij}(v)=0$.

Denote by $G(z,z')$ the Green's function of the Laplace operator $\Delta$ in $\End E_{0}$ ---
the integral kernel of the operator
$\Delta_{0}^{-1}$. It is an $\End\End\CC^{k}$-valued function on $\HH\times\HH$ with the transformation law
$$G(\s_{1} z,\s_{2} z')=\Ad\rho(\s_{1})G(z,z')\Ad\rho(\s_{2})^{-1},$$
where $\s_{1},\s_{2}\in\Gamma$ and $z,z'\in\HH$.
The Green's function is smooth when $z\neq\s z'$, $\s\in\Gamma$, and when $z'\rightarrow z$ it has a logarithmic singularity:
$$G(z,z')=-\frac{1}{\pi}\log|z-z'|\cdot I+O(1),\quad z'\rightarrow z,$$
where $I$ is now the identity element in $\End\End\CC^{k}$. The Green's function $Q(z,z')$ of the operator $\Delta$ in the trivial bundle
$\HH\times\End\CC^{k}$ is given by the explicit formula
$$Q(z,z')=-\frac{1}{\pi}\log\left|\frac{z-z'}{\z-z'}\right|\cdot I.$$
Set
\begin{equation} \label{psi}
\psi(z)=\frac{\pa}{\pa z'}\left.(G(z,z')-Q(z,z'))\right|_{z'=z}.
\end{equation}
The following results  will be used in Section \ref{LIT}.

\begin{lemma} \label{psi-expansion} 
The function $\psi:\HH\rightarrow\End\End\CC^{k}$ is smooth and  satisfies the transformation law
$$\psi(\s z)\s^{\prime}(z)=\Ad\rho(\s)\psi(z)\Ad\rho(\s)^{-1},\quad\s\in\Gamma;$$
in other words, $\psi\in C^{1,0}(X_{0},\End\End E_{0})$. Moreover, put
$$C_{i}=\lim_{y\rightarrow\infty}\psi(\sigma_{i}z)\sigma^{\prime}_{i}(z),\quad i=1,\dots,n.$$
Then we have 
$$C_{i}=(U_{i}\otimes \bar{U}_{i})T_{i}(U_{i}^{-1}\otimes \bar{U}_{i}^{-1}),$$
where $U_i$ is a unitary matrix that diagonalize $\rho(S_i)$, $\rho(S_{i})=U_{i}D_{i}U_{i}^{-1}$, and $T_{i}$ is the diagonal
$k^{2}\times k^{2}$-matrix with elements 
$$(T_i)_{lm}=-\mathrm{sgn}(\alpha^{i}_{l}-\alpha^{i}_{m})\sqrt{-1}\left(\frac{1}{2} -
|\alpha^{i}_{l}-\alpha^{i}_{m}|\right), \quad l,m=1,\dots, r_{i},$$ each repeated $k^{i}_{l}k^{i}_{m}$ times (we assume $\mathrm{sgn}(0)=1$).
\end{lemma}

\begin{proof}
The resolvent kernel of the operator $\Delta$  in the trivial bundle
$\HH\times\End\CC^{k}$, or, equivalently, the integral kernel of the operator $(\Delta +\frac{1}{2}s(s-1))^{-1}$, 
is $Q_{s}(z,z')=Q_{s}^{(0)}(z,z')\cdot I$, where $Q_{s}^{(0)}(z,z')$ is the resolvent kernel of 
$\Delta$ in the trivial line bundle
$\HH\times\CC$ (see, e.g., \cite{Takhtajan-Zograf91}\footnote{The Laplace operator in  \cite{Takhtajan-Zograf91} 
is $\frac{1}{2}\Delta$.}). In particular, 
$Q_{1}(z,z')=Q(z,z')$. The resolvent kernel $G_{s}(z,z')$ of the Laplace operator in the bundle $\End E^{\rho}$
is given by the series
$$G_{s}(z,z')=\sum_{\s\in\Gamma}Q_{s}(z,\s z')\Ad\rho(\s),$$  
that converges absolutely and uniformly on compact subsets of $\HH\times\HH$ for $\re s >1$.
We have 
$$G(z,z')=\lim_{s\rightarrow 1}\left(G_{s}(z,z')-\frac{1}{s(s-1)}\cdot\frac{1}{\pi k(2g-2+n)}\cdot I\right),$$
so that $\psi(z)$ given by \eqref{psi} is a smooth $\End\End E_{0}$-valued (1,0)-form on $X_0$,
or, equivalently, $\psi\in C^{1,0}(X_{0},\End\End E_{0})$. Now,
similar to Lemma 1 in \cite{Takhtajan-Zograf91} we have
$$\psi(\sigma_{i}z)\sigma^{\prime}_{i}(z)=\sum_{\substack{m=-\infty\\m\neq 0}}^{\infty}
\frac{\pa}{\pa z'}Q(z,z'+m)\Ad\rho(S_{i}^{m}) + o(1)\quad\text{as $y\rightarrow\infty$}.$$ 
Using the simple explicit formula
$$\frac{\pa}{\pa z'}Q(z,z')=\frac{1}{2\pi}\left(\frac{1}{z-z'} -\frac{1}{\z-z'}\right)\cdot I$$
we obtain that
$$C_{i}=\lim_{y\rightarrow\infty}\psi(\sigma_{i}z)\sigma^{\prime}_{i}(z)=\frac{1}{2\pi}\lim_{y\rightarrow\infty}\sum_{\substack{m=-\infty\\m\neq 0}}^{\infty}\left(\frac{1}{m+2\sqrt{-1}y}-\frac{1}{m}\right)\Ad\rho(S_{i}^{m}).$$
Choose a basis in $\CC^{k}$ such that $\rho(S_{i})$ is given by the diagonal matrix $D_{i}$. 
Using the elementary formulas
\begin{align*}
\sum_{m=1}^{\infty}\frac{2z}{z^{2}-m^{2}} & =\pi\cot\pi z-\frac{1}{z}, \quad z\in\CC\setminus \ZZ,\\
\intertext{and}
\sum_{m=1}^{\infty}\frac{\sin (2\pi m\alpha)}{m} &=\pi\left(\frac{1}{2}-\alpha\right),\quad 0<\alpha<1,
\end{align*}
with an odd extension to $-1<\alpha<0$, 
we easily get the second statement of the lemma.
\end{proof}

For each pair of harmonic forms $\mu,\nu\in\cH^{0,1}(X_{0},\End E_{0})$ we define a smooth $L^2$-section
$f_{\mu\bar\nu}\in\mathfrak{H}^{0,0}_{0}(X_{0},\mathrm{ad}\, E_{0})$ by the formula
$$f_{\mu\bar\nu}=\Delta_{0}^{-1}(\ast[\ast\mu,\nu]).$$

\begin{lemma}\label{f-mu-nu}
The section $f_{\mu\bar\nu}$ has the following asymptotics at the cusps:
$$f_{\mu\bar\nu}(\sigma_{i}z)=F_{\mu\bar\nu}^{i}+o(1)\quad \text{as}\quad y\rightarrow\infty,\quad i=1,\dots,n,$$
where $F_{\mu\bar\nu}^{i}\in\End\CC^{k}$ and $\tr F_{\mu\bar\nu}^{i}=0$. Moreover, for any $v\in V_{i}$ with $\tr v=0$ we have
\begin{align*}
\tr (F_{\mu\bar\nu}^{i}v)
& =2\int\limits_{X_{0}}
\ast[\ast\mu,\nu]\wedge \ast E_{i}(\,\cdot\,,v^{\ast};1) \\
& = 4\iint\limits_{F}\tr([\mu(z),\nu(z)^{\ast}]E_{i}(z,v,1))dxdy,\\
\end{align*}
whereas $\tr (F_{\mu\bar\nu}^{i}v)=0$ for $v\notin V_{i}$.
\end{lemma}
\begin{proof}
We repeat the main steps of the proof of Lemma 2 in Section 1 of \cite{Takhtajan-Zograf91}.
%\footnote{The second 
%line of the first formula on p.~408 in \cite{Takhtajan-Zograf91} misses the overall negative sign.}. 
We have $\Ad\rho(S_{i})v=e^{2\sqrt{-1}\pi\beta}v$ with some $\beta\in (-1,1)$. Put $g(z)=\tr(f_{\mu\bar\nu}(\sigma_{i}z)v)$. 
The function $g(z)$ has the property $g(z+1)=e^{2\sqrt{-1}\pi\beta}g(z)$, so it admits the Fourier expansion 
$$g(z)=\sum_{m=-\infty}^{\infty}a_{m}(y)e^{2\pi\sqrt{-1}(m+\beta)x},\quad z\in\HH.$$ 
Automorphic forms
$\mu(z), \nu(z)$ are exponentially decreasing at the cusps, and the function $g(z)$ is square integrable on $F$ 
with respect to the hyperbolic area form. From the equation  
$\Delta g=\tr(\ast[\ast\mu,\nu]v)$ it then follows that the functions $\displaystyle{\frac{d^{2}a_{m}}{dy^{2}}-4\pi^{2}(m+\beta)^{2}}$ for all $m\in\ZZ$ are exponentially decreasing as $y\rightarrow\infty$. Thus, when $\beta\neq 0$, the function $g(z)$ decays exponentially as $y\rightarrow\infty$. 
When $\beta=0$, the coefficient
$a_{0}(y)=a_{0}$ is a constant, and
$g(z)-a_{0}$ exponentially decays as $y\rightarrow\infty$. 
To get the integral formula for $a_{0}$, consider 
a canonical fundamental domain $F$ for $\Gamma$ with exactly $n$ cusps at the points $x_{1},\dots,x_{n}$, 
and take $F^{Y}=\{z\in F | \im(\sigma_{i}^{-1}z)\leq Y,\, i=1,\dots,n\}$. Using Green's formula 
and asymptotics \eqref{E-asymptotic}, we get
\begin{gather*}
2\iint\limits_{F}\tr([\mu(z),\nu(z)^{\ast}]E_{i}(z,v,1))dxdy =\iint\limits_{F}\tr(\Delta_{0}f_{\mu\bar\nu}(z)E_{i}(z,v,1))\frac{dxdy}{y^{2}} \\ 
 =\lim_{Y\rightarrow\infty}\frac{1}{2}\int\limits_{\pa F^{Y}}\tr\left\{
E_{i}(z,v;1)\left(\frac{\pa f_{\mu\bar\nu}}{\pa y}dx -   \frac{\pa f_{\mu\bar\nu}}{\pa x}dy \right)\right. \\
\left.-f_{\mu\bar\nu}\left( \frac{\pa}{\pa y}E_{i}(z,v;1)dx-\frac{\pa}{\pa x}E_{i}(z,v;1)dy\right)\right\}=\frac{a_{0}}{2},
\end{gather*}
where we used the differential equation \eqref{E-Laplace} for $s=1$. 
Note that by definition $a_{0}=\tr (F_{\mu\bar\nu}^{i}v)$, which completes the proof.
\end{proof}

\section{The moduli space of parabolic bundles} \label{moduli space}

\subsection{The complex structure}\label{cs}

According to the Mehta-Seshadri theorem, the moduli space $\mathcal{N}$ of stable parabolic bundles of rank $k$ 
on $X=\Gamma\bk\overline\HH$ with given weights and multiplicities at the marked points $P_1,\dots,P_n\in X$
is isomorphic to the space $\Hom(\Gamma,U(k))^{0}/U(k)$ of equivalence classes of irreducible admissible representations of $\Gamma$ (where the unitary group $U(k)$ acts by conjugation).   
This is a complex manifold of dimension   
$$d=k^{2}(g-1) +1+\sum_{i=1}^{n}\dim_{\CC}\cF_{i}.$$
If the parabolic structure is \emph{integral} (i.e.
$\sum_{l=1}^{r(P)}k_{l}(P)\alpha_{l}(P)\in\ZZ$ for each $P\in S$) one can consider unimodular irreducible admissible
representations of $\Gamma$. The representation space $\mathcal{N}_{0}=\Hom(\Gamma,SU(k))^{0}/SU(k)$ is then a complex
submanifold of $\mathcal{N}$ of dimension $d_{0}=d-g$. 
The correspondence
$E\mapsto\wedge^{k}E$ defines a holomorphic mapping $\mathcal{N}\rightarrow J_{\mathrm{deg} E}$, 
where $J_{\mathrm{deg} E}$ is the component of the Picard group $\mathrm{Pic}(X)$ parametrizing line bundles  
of degree $\mathrm{deg}\,E$ on $X$. The fibers of this mapping --- the moduli spaces of stable parabolic vector bundles on $X$ with fixed determinant --- are all isomorphic to $\mathcal{N}_{0}$ as complex manifolds.

As in \cite{Narasimhan-Seshadri64}, the holomorphic tangent space $T_{\{E\}}\mathcal{N}$ 
at the point $\{E\}\in\mathcal{N}$ corresponding to the stable parabolic bundle $E$ is identified with the space $\cH^{0,1}(X_{0},\End E_{0})$ of square integrable harmonic
$(0,1)$-forms on $X_{0}$ with values in $\End E_{0}$. The corresponding holomorphic cotangent space $T^{\ast}_{\{E\}}\mathcal{N}$ is identified with the space  $\cH^{1,0}(X_{0},\End E_{0})$ of square integrable harmonic $(1,0)$-forms on $X_{0}$ with values in $\End E_{0}$, and the pairing 
$$\cH^{0,1}(X_{0},\End E_{0})\otimes\cH^{1,0}(X_{0},\End E_{0})\rightarrow\CC$$ 
is given by
$$(\nu,\theta)\mapsto\int\limits_{X_{0}}\nu\wedge\theta, \quad\nu
\in\cH^{0,1}(X_{0},\End E_{0}),\,\theta\in\cH^{1,0}(X_{0},\End E_{0}).$$

Let $\rho:\Gamma\rightarrow U(k)$ be an admissible irreducible representation. Exactly as in \cite{Takhtajan-Zograf89}, 
we can show that, for each $\nu\in\cH^{0,1}(X_{0},\End E_{0})$ sufficiently close to zero, there exists a unique mapping
$f^{\nu}: \HH\rightarrow\mathrm{GL}(k,\CC)$ with the following properties:
\begin{enumerate}
\item[(i)] $f^\nu$ satisfies the equation $$\frac{\pa f^{\nu}}{\pa \z}=f^{\nu}(z)\nu(z),\quad z\in\HH;$$
\item[(ii)] $\det f^{\nu}(z_{0})=1$ at some fixed $z_{0}\in\HH$
(say, $z_{0}=\sqrt{-1}$);
\item[(iii)] $\rho^{\nu}(\gamma)=f^{\nu}(\gamma z)\rho(\gamma)f^{\nu}(z)^{-1}$
is independent of $z$ and is an admissible irreducible unitary representation of $\Gamma$; 
\item[(iv)] $f^{\nu}$ is regular at the cusps, that is,
$$f^{\nu}(x_{i})=\lim_{z\rightarrow\infty}f^{\nu}(\sigma_{i}z)<\infty,\quad i=1,\dots,n.$$
\end{enumerate}

Let $\nu_{1},\dots,\nu_{d}$ be a basis for $\cH^{0,1}(X_{0},\End E_{0}^{\rho})$, and let $\nu=\vep_{1}\nu_{1}+\dots+\vep_{d}\nu_{d}$, where
$\vep_{i}\in\CC,\; i=1,\dots,d,$ are sufficiently small.
The mapping $(\vep_{1},\dots,\vep_{d})\mapsto \{E^{\rho^{\nu}}\}$ provides a coordinate chart on $\mathcal{N}$ in the neighborhood of the point $\{E^{\rho}\}$.
These coordinates transform holomorphically and endow $\mathcal{N}$ with the structure of a complex manifold (they are 
similar to Bers' coordinates on Teichm\"{u}ller spaces). The differential of such coordinate transformation is a linear mapping $\cH^{0,1}(X_{0},\End E_{0}^{\rho})\rightarrow\cH^{0,1}(X_{0},\End E_{0}^{\rho^{\nu}})$ explicitly given by the formula
\begin{equation} \label{differential}
\mu\mapsto P_{\nu}(\Ad f^{\nu}(\mu)),\quad\mu\in\cH^{0,1}(X_{0},\End E_{0}^{\rho}).
\end{equation}
Here $P_{\nu}: \mathfrak{H}^{0,1}(X_{0},\End E_{0}^{\rho^{\nu}})\rightarrow\cH^{0,1}(X_{0},\End E_{0}^{\rho^{\nu}})$ is the orthogonal projection, and $\Ad f^{\nu}$ is understood as a fiberwise linear mapping $\End E_{0}^{\rho}\rightarrow\End E_{0}^{\rho^{\nu}}$, where $\Ad f^{\nu}(\mu)=f^{\nu}\cdot\mu\cdot (f^{\nu})^{-1}$. When the parabolic structure is integral, 
the holomorphic tangent space $T_{\{E^{\rho}\}}\mathcal{N}_{0}$
at the point $\{E^{\rho}\}\in\mathcal{N}_{0}$ is identified with the subspace 
$\cH^{0,1}(X_{0},\ad E_{0}^{\rho})\hookrightarrow\cH^{0,1}(X_{0},\End E_{0}^{\rho})$. 
Note that there is an orthogonal decomposition
$$\cH^{0,1}(X_{0},\End E_{0}^{\rho})\cong\cH^{0,1}(X_{0},\ad E_{0}^{\rho})\oplus\cH^{0,1}(X_{0})\otimes I.$$
If the basis $\nu_{1},\dots,\nu_{d}$ for $\cH^{0,1}(X_{0},\End E_{0}^{\rho})$
is chosen in such a way that
$\nu_{1},\dots,\nu_{d_{0}}\in\cH^{0,1}(X_{0},\ad E_{0}^{\rho})$ and $\nu_{d_{0}+1},\dots,\nu_{d}\in\cH^{0,1}(X_{0})\otimes I$, then in the local coordinates 
$(\vep_{1},\dots,\vep_{d})$ the submanifold $\mathcal{N}_{0}\subset\mathcal{N}$ is given by the equations $\vep_{d_{0}+1}=\cdots=\vep_{d}=0$.

The moduli space $\mathcal{N}$ carries a Hermitian metric given by the inner product \eqref{inner-1} in the fibers of $T\mathcal{N}$.
This metric is analogous to the Weil-Petersson metric on Teichm\"{u}ler space, and for the moduli spaces of stable bundles of fixed rank and degree was introduced in \cite{Narasimhan69, Atiyah-Bott83}. This metric is K\"{a}hler and we will denote its K\"{a}hler (symplectic) form by $\Omega$:
$$\Omega\left(\frac{\pa}{\pa \vep(\mu)},\frac{\pa}{\pa\overline{ \vep(\nu)}}\right)=\frac{\sqrt{-1}}{2}\langle\mu,\nu\rangle.$$
Here $\displaystyle{\frac{\pa}{\pa \vep(\mu)}}$ and $\displaystyle{\frac{\pa}{\pa\overline{ \vep(\nu)}}}$ are the holomorphic and antiholomorphic tangent vectors at  $\{E\}\in\mathcal{N}$ corresponding to $\mu,\nu\in\cH^{0,1}(X_{0},\End E_{0})$ respectively.

\subsection{Families of endomorphism bundles}

It follows from the general deformation theory that the moduli space $\mathcal{N}$ admits an open covering $\mathcal{N}=\cup_{\alpha\in A}U_{\alpha}$
such that for every $\alpha\in A$ there exists a \emph{family of endomorphism bundles} on $X_{0}\times U_{\alpha}$: 
a holomorphic vector bundle $\cE_{\alpha}\rightarrow X_{0}\times U_{\alpha}$ with the Hermitian metric $h_{\cE_{\alpha}}$
such that $\left.\cE_{\alpha}\right|_{X_{0}\times\{E\}}\cong\End E_{0}$ as Hermitian vector bundles for any $\{E\}\in U_{\alpha}$. If we consider only traceless endomorphisms, we get a family of the adjoint bundles $\cF_{\alpha}$ for which $\left.\cF_{\alpha}\right|_{X_{0}\times\{E\}}\cong\ad E_{0}$ for $\{E\}\in U_{\alpha}$, and $\cE_{\alpha}=\cF_{\alpha}\oplus\CC$, where $\CC$ is understood as the trivial line bundle on $X_{0}\times\ U_{\alpha}$.

The direct image $\pi_{\ast}\cE_{\alpha}$ of $\cE_{\alpha}$ under the projection $\pi: X_{0}\times U_{\alpha}\rightarrow\
U_{\alpha}$ is isomorphic to the restriction $\left.T\mathcal{N}\right|_{U_{\alpha}}$ of the tangent bundle $T\mathcal{N}$ 
to $U_{\alpha}$. Correspondingly, $\left.T^{\ast}\mathcal{N}\right|_{U_{\alpha}}\cong \pi_{\ast}(\left.\cE_\alpha\otimes T^{\ast}_{V}\right|_{X_{0}\times U_{\alpha}})$, where $T^{\ast}_{V}$ is the vertical
(along the fibers of the projection $\pi$) cotangent bundle on $X_{0}\times U_{\alpha}$. If an open covering is chosen
properly, then for every $U_{\alpha}$ there exist $d$
holomorphic sections $\omega_{1},\dots,\omega_{d}$ of $\cE_{\alpha}\otimes T_{V}^{\ast}$ on $X_{0}\times U_{\alpha}$ that are linearly independent over each fiber $X_{0}\times\{E\}$, $\{E\}\in U_{\alpha}$. This means, in particular, that over each
point $\{E\}=\{E^{\rho}\}\in U_{\alpha}$ the sections $\left.\omega_{1}\right|_{X_{0}\times\{E^{\rho}\}},\dots,\left.\omega_{d}\right|_{X_{0}\times\{E^{\rho}\}}$ of $\End E_{0}^{\rho}\otimes T^{\ast}X_{0}$ form a basis for the vector space $\cH^{1,0}(X_{0},\End E_{0}^{\rho})$ and for every $\nu\in\cH^{0,1}(X_{0},\End E_{0}^{\rho})$ each of the forms
$$\Ad(f^{\vep\nu})^{-1}(\left.\omega_{i}\right|_{X_{0}\times\{E^{\rho}\}})\in C^{1,0}(X_{0},\End E_{0}^{\rho}),\quad i=1,\dots,d,$$
is holomorphic in $\vep\in\CC$ at $\vep=0$. 

For the integral parabolic structure put $V_{\alpha}=U_{\alpha}\cap\mathcal{N}_{0}$. Then we have 
$\pi_{\ast}(\left.\cF_{\alpha}\right|_{X_{0}\times V_{\alpha}})\cong \left.T\mathcal{N}_{0}\right|_{V_{\alpha}}$  and $\pi_{\ast}(\left.\cF_{\alpha}\otimes T^{\ast}_{V}\right|_{X_{0}\times V_{\alpha}})\cong \left.T^{\ast}\mathcal{N}_{0}\right|_{V_{\alpha}}$.
The sections $\omega_{1},\dots,\omega_{d}$ of the bundle $\cE_{\alpha}\otimes T_{V}^{\ast}$ can be chosen in such 
a way that $\omega_{1},\dots,\omega_{d_{0}}$ take values in the subbundle 
$\cF_{\alpha}\otimes T_{V}^{\ast}\hookrightarrow\cE_{\alpha}\otimes T_{V}^{\ast}$.

\begin{remark}
To the best of our knowledge, it is not completely clear whether there 
always exists a \emph{universal endomorphism bundle} $\cE\rightarrow X_{0}\times\mathcal{N}$
such that $\left.\cE\right|_{X_{0}\times\{E\}}\cong \End E_{0}$ for every $\{E\}\in\mathcal{N}$.
For a generic weight system the existence of the universal endomorphism bundle follows e.g. from
\cite{BY}, Proposition 3.2.
\end{remark}

\section{Variational formulas}

\subsection{Lie derivatives} \label{lie}

By definition, a \emph{family of forms of type} $(p,q)$, $p,q=0,1$, on $X_{0}\times U_{\alpha}$ is a smooth section of the bundle $\cE_{\alpha}\otimes\wedge^{p}T_{V}^{\ast}\otimes\wedge^{q}\overline{T}_{V}^{\ast}\rightarrow X_{0}\times U_{\alpha}$, where $T_{V}^{\ast}$ and $\overline{T}_{V}^{\ast}$ are the holomorphic and antiholomorphic vertical cotangent bundles 
on $X_{0}\times U_{\alpha}$ respectively. Let $\{E^{\rho^{\vep\nu}}\}$ for sufficiently small $\vep\in\CC$ be a complex 
curve in $\mathcal{N}$ with the tangent vector $\pa/\pa\vep(\nu)$ at the point $\{E^{\rho}\}\in U_{\alpha}$, where $\nu\in\cH^{0,1}(X_{0},\End E_{0}^{\rho})$, and let $\omega^{\vep}\in C^{p,q}(X_{0},\End E_{0}^{\rho^{\vep\nu}})$ 
be a family of forms of type $(p,q)$ over this curve. The Lie derivatives of the family
$\omega^{\vep}$ in the directions $\pa/\pa\vep(\nu)$ and $\pa/\pa\overline{\vep(\nu)}$ are defined by the standard formulas
$$L_{\nu}\omega =\left.\frac{\pa}{\pa\vep}\right|_{\vep=0} \Ad (f^{\vep\nu})^{-1}(\omega^{\vep}),\quad L_{\bar{\nu}}\omega =\left.\frac{\pa}{\pa\bar{\vep}}\right|_{\vep=0} \Ad (f^{\vep\nu})^{-1}(\omega^{\vep}).$$
The Lie derivatives of smooth families of linear operators 
$$A^{\vep}:
\mathfrak{H}^{p,q}(X_{0},\End E_{0}^{\rho^{\vep\nu}})\rightarrow\mathfrak{H}^{p',q'}(X_{0},\End E_{0}^{\rho^{\vep\nu}})$$ 
are defined by the formulas
\begin{align*}
L_{\nu}A &=\left.\frac{\pa}{\pa\vep}\right|_{\vep=0} \Ad (f^{\vep\nu})^{-1} \circ A^{\vep} \circ\Ad f^{\vep\nu},\\
L_{\bar{\nu}}A &=\left.\frac{\pa}{\pa\bar{\vep}}\right|_{\vep=0} \Ad (f^{\vep\nu})^{-1} \circ A^{\vep}\circ \Ad f^{\vep\nu}.
\end{align*}
These are linear operators from $\mathfrak{H}^{p,q}(X_{0},\End E_{0}^{\rho})$ to $\mathfrak{H}^{p',q'}(X_{0},\End E_{0}^{\rho})$. The Lie derivatives obey the Leibniz rules; in particular,
$$L_{\nu}(A\omega)=(L_{\nu}A)\omega +A(L_{\nu}\omega),\quad L_{\bar\nu}(A\omega)=(L_{\bar\nu}A)\omega +A(L_{\bar\nu}\omega).$$

Repeating verbatim the computations in \cite{Takhtajan-Zograf89}, we get the formulas
\begin{equation} \label{metric-first}
L_{\nu}h_{\cE_{\alpha}}(\xi,\eta)=L_{\bar\nu}h_{\cE_{\alpha}}(\xi,\eta)=0,
\end{equation}
\begin{equation} \label{metric-second}
L_{\mu}L_{\bar\nu}h_{\cE_{\alpha}}(\xi,\eta)=-h_{\cE_{\alpha}}([f_{\mu\bar{\nu}},\xi],\eta)
\end{equation}
for all $\mu,\nu\in\cH^{0,1}(X_{0},\End E_{0}^{\rho})$ and all bounded $\xi,\eta\in C^{0}(X_0, \End E_0^{\rho})$;
here $f_{\mu\bar{\nu}}=\Delta^{-1}_{0}(\ast[\ast\mu,\nu])$ as in Lemma \ref{f-mu-nu}. Furthermore,
\begin{alignat*}{2}
L_{\nu}\delb &=\ad\nu,&\quad L_{\bar\nu}\delb & =0,\\
L_{\nu}\delb^{\ast}&=0,&\quad L_{\bar\nu}\delb^{\ast} & =-\ast\ad\ast\nu.
\end{alignat*}
so that for the operators $\Delta=\delb^{\ast}\delb$ and $P=I-\delb\Delta_{0}^{-1}\delb^{\ast}$ we get
\begin{equation}\label{var-delta}
L_{\nu}\Delta=\delb^{\ast}\ad\nu\quad\text{and}\quad L_{\bar\nu}P=\delb\Delta_{0}^{-1}\ast\ad\nu\ast P.
\end{equation}
For the family $\mu^{\vep\nu}=P_{\vep\nu}(\Ad f^{\vep\nu}\mu)$ which corresponds, under the identification
$T_{\{E^{\rho^{\vep\nu}}\}}\mathcal{N}\cong\cH^{0,1}(X_{0},\End E_{0}^{\rho^{\vep\nu}})$, to the tangent 
vector field $\pa/\pa\vep(\mu)$ to the (complex) curve $\{E^{\rho^{\vep\nu}}\}\in\mathcal{N}$, 
we get
\begin{equation} \label{var-mu}
L_{\bar\nu}\mu=\delb f_{\mu\bar\nu}.
\end{equation}

The determinant of the K\"{a}hler metric on $\mathcal{N}$ is a Hermitian metric in the canonical line bundle $\det T^{\ast}\mathcal{N}=\wedge^{d}T^{\ast}\mathcal{N}$. Its curvature $(1,1)$-form $\Theta$ is given by
\begin{equation}
\Theta\!\left(\frac{\pa}{\pa \vep(\mu)},\frac{\pa}{\pa\overline{ \vep(\nu)}}\right)=-\Tr((\ad f_{\mu\bar\nu}I + (L_{\mu}\delb)\Delta_{0}^{-1}(L_{\bar\nu}\delb^{\ast}))P),
\end{equation}  
where $\Tr$ is the operator trace in the Hilbert space $\mathfrak{H}^{0,1}(X_{0}, \End E_{0})$, $\ad f_{\mu\bar\nu}$
is a linear operator in $C^{0,1}(X_{0}, \End E_{0})$ understood as
$$\ad f_{\mu\bar\nu}(\xi)=[f_{\mu\bar\nu},\xi],\quad \xi\in C^{0,1}(X_{0}, \End E_{0}),$$ 
$I$ is the identity operator in $\mathfrak{H}^{0,1}(X_{0},\End E_{0})$, and $P:\mathfrak{H}^{0,1}(X_{0},\End E_{0})\rightarrow \cH^{0,1}(X_{0},\End E_{0})$ is the orthogonal projection.

\subsection{Eisenstein-Maass series and closed $(1,1)$-forms}

Let $E\cong E^{\rho}$ be a stable parabolic vector bundle on $X$.
For each marked point $P_{i}\in X$ and each vector $v\in V_{i}=\ker(\Ad\rho(S_{i})-I)$ with $\tr v=0$, $i=1,\dots,n$, 
we define a $(1,1)$-form $\Omega_{i,v}$ in a neighborhood of the point $\{E\}\in \mathcal{N}$ as follows. Choose a 
basis $\varphi_1,\dots,\varphi_d\in\cH^{0,1}(X_{0},\End E_{0}^{\rho})$ and put $\varphi=\vep_1\varphi_1+\dots +\vep_d\varphi_d$
with small enough $\vep_1,\dots,\vep_d\in\CC$. The parameters $\vep_1,\dots,\vep_d$ provide local
coordinates near the point $\{E\}\in \mathcal{N}$ by means of the mapping 
$$(\vep_1,\dots,\vep_d)\mapsto\{E^{\rho^\varphi}\}\in\mathcal{N}$$ 
(see Section \ref{cs} for details).

For any $\mu,\nu\in\cH^{0,1}(X_{0},\End E_{0}^{\rho})$ consider two families of harmonic 
(0,1)-forms $\mu^\varphi,\,\nu^\varphi\in\cH^{0,1}(X_{0},\End E_{0}^{\rho^\varphi})$,
where $\mu^\varphi=P_\varphi(\Ad f^\varphi(\mu))$ and $\nu^\varphi=P_\varphi(\Ad f^\varphi(\nu))$.
At the point $\{E\}\in \mathcal{N}$ the form $\Omega_{i,v}$ is defined by the formula
\begin{align*}
\Omega_{i,v}\left(\frac{\pa}{\pa \vep(\mu)},\frac{\pa}{\pa\overline{ \vep(\nu)}}\right) & =\frac{\sqrt{-1}}{2}\langle\ast[\ast\mu,\nu],E_{i}(\,\cdot\,,v^{\ast};1)\rangle \\
& =\sqrt{-1}\iint\limits_{F}\tr\left([\mu(z),\nu(z)^{\ast}]E_{i}(z,v;1)\right)dxdy.
\end{align*}
It extends to the neighborhood of $\{E\}\in\mathcal{N}$ by replacing $\mu,\nu$ with $\mu^{\varphi},\nu^{\varphi}$,
and $v\in V_{i}$ with $v^{\varphi}=f^{\varphi}v\in V_{i}^{\varphi}=\ker(\Ad\rho^{\varphi}(S_{i})-I)$. 
Note that by Lemma \ref{f-mu-nu} we also have
\begin{equation} \label{Omega-F}
\Omega_{i,v}\left(\frac{\pa}{\pa \vep(\mu)},\frac{\pa}{\pa\overline{ \vep(\nu)}}\right)=\frac{\sqrt{-1}}{4}\tr(F^{i}_{\mu\bar\nu}v).
\end{equation}

\begin{lemma} \label{Closed}
The $(1,1)$-forms $\Omega_{i,v}$ on $\mathcal{N}$ are closed and satisfy the condition $\overline{\Omega}_{i,v}=\Omega_{i,v^{\ast}}$.
\end{lemma}

\begin{proof} To get the equality $d\Omega_{i,v}=0$ it is sufficient to show that 
$$\frac{\pa}{\pa\vep(\mu)}\langle\ast[\ast\nu,\lambda],E_{i}(\,\cdot\,,v;1)\rangle = \frac{\pa}{\pa\vep(\nu)}\langle\ast[\ast\mu,\lambda],E_{i}(\,\cdot\,,v;1)\rangle$$
for all  $\mu,\nu,\lambda\in\cH^{0,1}(X_{0},\End E_{0})$. It can be verified exactly as in Lemma 3 
of \cite{Takhtajan-Zograf91} using formulas \eqref{var-delta}, \eqref{var-mu} and the equality
$$L_{\mu}E_{i}(\,\cdot\,,v;1)=\Delta_{0}^{-1}(\delb^{\ast}\ad\mu\, E_{i}(\,\cdot\,,v;1)),$$
which follows from \eqref{var-delta}. To verify the complex conjugation property, we observe that, since
$f_{\mu\bar\nu}(z)^{\ast}=f_{\nu\bar\mu}(z)$, we have $(F^{i}_{\mu\bar\nu})^{\ast}=F^{i}_{\nu\bar\mu}$, and from the cyclic
invariance of the trace we get
$\overline{\tr(F^{i}_{\mu\bar\nu}v)}=\tr(F^{i}_{\mu\bar\nu}v)^{\ast}=
\tr(F^{i}_{\nu\bar\mu}v^{\ast})$.
\end{proof}

Let $u_{j}$, $j=1,\dots, k_{l}=k^{i}_{l}$, be an orthonormal basis for the eigenspace of $\rho(S_{i})$ in $\CC^{k}$ corresponding to the eigenvalue $e^{ 2\pi\sqrt{-1}\alpha^{i}_{l}}$, and put 
$$v_{j}=u_{j}\otimes \bar{u}_{j}-I/k\in\End\CC^{k},\quad \tr v_{j}=0.$$
Since $v_{j}^{\ast}=v_{j}$, the $(1,1)$-forms $\Omega_{i,v_{j}}$ are real. 
Moreover, because $\sum_{j=1}^{k_{l}}u_{j}\otimes \bar{u}_{j}$ represents an orthogonal projection to the eigenspace
corresponding to the eigenvalue $e^{ 2\pi\sqrt{-1}\alpha^{i}_{l}}$ of $\rho(S_{i})$, the $(1,1)$-forms
$$\Omega_{il}=\sum_{j=1}^{k_{l}}\Omega_{i,v_{j}}$$
do not depend on the choice of the basis $\{u_{j}\}_{j=1}^{k_{l}}$ and are well-defined on the moduli space $\mathcal{N}$.

\subsection{Holomorphic line bundles} \label{bundles}

Here we realize closed, real $(1,1)$-forms $\Omega_{il}$ as the curvature forms (more precisely, as the first Chern forms) 
of certain natural line bundles on the moduli space $\mathcal{N}$. Namely, for each $P_{i}\in S$ and $l=1,\dots,r_{i}$, let
$\lambda_{il}$ be the holomorphic line bundle on $\mathcal{N}$ whose fiber over the point $\{E\}\in\mathcal{N}$ is the complex
line $\det W_{il}$, where $W_{il}=F_{l}E_{P_{i}}/F_{l+1}E_{P_{i}}$ is the complex vector space of dimension 
$\dim W_{il}=k_{l}=k_{l}^{i}$. We introduce a Hermitian metric $\Vert\cdot\Vert_{il}$ in the line bundle $\lambda_{il}$ as
follows. By the Mehta-Seshadri theorem we have $E\cong E^{\rho}$, where $\rho$ is an irreducible admissible representation 
of the group $\Gamma$. As in the previous section, let $u_{1},\dots,u_{k_{l}}$ be orthonormal eigenvectors of the unitary matrix $\rho(S_{i})$ corresponding to the eigenvalue $e^{2\pi\sqrt{-1}\alpha_{l}^{i}}$. Then the Hermitian metric $\Vert\cdot\Vert_{il}$ is defined by the standard Hermitian norm of the vector $u=u_{1}\wedge\dots\wedge u_{k_{l}}\in\wedge^{k_{l}}\CC^{k}$,
$$\Vert u\Vert_{il}^{2}=\det\{(u_{j},u_{l})\}_{j,m=1}^{k_{l}}=1.$$

\begin{lemma} Let $c_{1}(\lambda_{il}, \Vert\cdot\Vert_{il})$ denote the first Chern form of the line bundle $\lambda_{il}$ with respect to the metric $\Vert\cdot\Vert_{il}$. Then
$$c_{1}(\lambda_{il}, \Vert\cdot\Vert_{il})=\frac{2}{\pi}\,\Omega_{il},\quad i=1,\dots,n,\;l=1,\dots,r_{i}.$$
\end{lemma}

\begin{proof}
For $\mu,\nu\in\curly{H}^{0,1}(X_{0},\mathrm{End}\, E_{0})$ 
let $\varphi=\vep_{1}\mu+\vep_{2}\nu$, and put
$$\Phi_{\mu\nu}(z)=\left.\frac{\partial^{2}}{\partial\vep_{1}\partial\bar{\vep}_{2}}\right|_{\vep_{1}=\vep_{2}=0}\left(f^{\varphi}(z)^{*}f^{\varphi}(z)\right)\in C^{0}(X_{0},\End E_{0}).$$
As in \cite{Takhtajan-Zograf89}, we obtain
$$\Delta\Phi_{\mu\nu}=-\ast[\ast\mu,\nu],$$
so that $\Phi_{\mu\nu}=-f_{\mu\nu} +c I$. Normalizing the mapping $f^{\varphi}$ as in Section \ref{lie}, we get
$\tr\Phi_{\mu\nu}(z_{0})=0$, so that $c=0$. 
Put
$$u_{j}^{\varphi}=\lim_{z\rightarrow x_{i}}f^{\varphi}(z)u_{j} =f^{\varphi}(x_{i})u_{j},\quad j=1,\dots, k_{l}.$$ 
Now for
$u^{\varphi}=u_{1}^{\varphi}\wedge\dots\wedge u_{k_{l}}^{\varphi}$ we have
$$\Vert u^{\varphi}\Vert^{2}_{il}= \det\{(f^{\varphi}(x_{i})^{*}f^{\varphi}(x_{i})u_{j},u_{l})\}_{j,m=1}^{k_{l}},$$ 
and using the fact that $\mathrm{tr}F^{i}_{\mu\bar\nu}=0$ we derive
\begin{align*}
\left.\frac{\partial^{2}}{\partial\vep_{1}\partial\bar{\vep}_{2}}\right|_{\vep_{1}=\vep_{2}=0}\log\Vert u^{\varphi}\Vert^{2}_{il}
& =-\lim_{z\rightarrow x_{i}}\sum_{j=1}^{k_{l}}(f_{\mu\bar\nu}(z)u_{j},u_{j}) \\
& =-\sum_{j=1}^{k_{l}}\mathrm{tr}(F^{i}_{\mu\bar\nu}v_{j}),
\end{align*}
so that the desired statement follows now from \eqref{Omega-F}.
\end{proof}

\begin{remark}
For the moduli space of punctured Riemann surfaces similar results were obtained
in \cite{Weng,Wolpert}. Here we use the approach of \cite{Wolpert}.
\end{remark}

\section{Local index theorems} \label{LIT}

\subsection{The first variation of the Selberg zeta function}

Recall (see \cite[Ch.~5]{Venkov81}) that the Selberg zeta function $Z(s,\Gamma;\chi)$ for the Fuchsian group $\Gamma$ with the unitary representation $\chi$ is defined for $\re s>1$ as the following absolutely convergent product
$$Z(s,\Gamma;\chi)=\prod_{\{\gamma\}}\prod_{k=0}^{\infty}\det(I-\chi(\s)N(\s)^{-s-k}),$$
where $\{\s\}$ runs over the set of all primitive conjugacy 
classes of hyperbolic elements of $\Gamma$, and $N(\s)>1$ is the norm of the element $\s\in\Gamma$, 
i.e., $\s$ is conjugate to the diagonal matrix $\left(\begin{smallmatrix}
N(\s)^{1/2} 
& 0\\ 0 & \;\;N(\s)^{-1/2}\end{smallmatrix}\right)$ in $\mathrm{SL}(2,\RR)$. 
The logarithmic derivative of the Selberg zeta function for $\re s>1$ is given by the integral
\begin{equation} \label{STF}
\frac{1}{2s-1}\frac{d}{ds}\log Z(s,\Gamma;\chi)=\frac{1}{2}\iint\limits_{F}\sum_{\text{$\s$ hyperbolic}}
\tr\chi(\s)Q_{s}(z,\s z)\frac{dxdy}{y^{2}},
\end{equation}
where the sum is taken over all hyperbolic elements in $\Gamma$
(see, e.g., \cite[Theorem 4.3.4, part 2)]{Venkov81}). The function $Z(s,\Gamma;\chi)$ is positive for 
$s\in (1,\infty)$ and admits a meromorphic continuation to the whole $s$-plane.

If $\chi=\Ad\rho$, where $\rho$ is an admissible irreducible representation, then $Z(s,\Gamma; \Ad\rho)$ has a simple 
zero at $s=1$, and as in \cite{Takhtajan-Zograf91} we define the regularized determinant of the Laplace operator 
$\Delta$ in $\mathfrak{H}^{0,0}(X_{0},\End E_{0}^{\rho})$ by the formula
$$\det\Delta=\left.\frac{\pa}{\pa s} \right|_{s=1}Z(s,\Gamma,\Ad\rho)
=\lim_{s\rightarrow\ 1} \frac{1}{s-1} Z(s,\Gamma,\Ad\rho).$$
As a function of $\rho\in\Hom(\Gamma,U(k))^{0}/U(k)\cong\mathcal{N}$, the determinant $\det\Delta$ is smooth
and positive. Denote by $\del_{\mathcal{N}}$ and $\delb_{\mathcal{N}}$ the $(1,0)$- and $(0,1)$-components of 
the de Rham differential on $\mathcal{N}$ respectively.
\begin{lemma} \label{1st-variation} Let $\mu\in\cH^{0,1}(X_{0},\End E_{0})$. Then at the point $\{E\}\in\mathcal{N}$
$$\del_{\mathcal{N}}\log\det\Delta\left(\frac{\pa}{\pa\vep(\mu)}\right)=-\sqrt{-1}\int\limits_{X_{0}}\ad\mu\wedge\psi,$$
where $\ad\mu=[\mu,\,\cdot\,]$ is understood as an element in $\cH^{0,1}(X_{0},\End\End E_{0})$, and $\psi\in C^{1,0}(X_{0},\End\End E_{0})$ is given by \eqref{psi}.
\end{lemma}
\begin{proof} We prove this lemma starting with the formula
$$\frac{\pa}{\pa\vep(\mu)}\log\det\Delta=\lim_{s\rightarrow 1^{+}}L_{\mu}\log Z(s,\Gamma;\Ad\rho),$$
and repeating the proof of Theorem 1 in \cite{Takhtajan-Zograf88} (see also Lemma 3 in Section 3 in \cite{Takhtajan-Zograf91})
with obvious adjustments for $Z(s,\Gamma;\Ad\rho)$.  
\end{proof}

\subsection{Local index theorem in Quillen's form}

Let $\tilde{\Omega}$ be the $(1,1)$-form on $\mathcal{N}$ defined
at each point $\{E\}\in\mathcal{N}$ by
\begin{equation} \label{omega}
\tilde\Omega\left(\frac{\pa}{\pa \vep(\mu)},\frac{\pa}{\pa\overline{ \vep(\nu)}}\right)=\frac{\sqrt{-1}}{2}\int\limits_{X_{0}}\ad\mu\wedge\ad\ast\nu,
\end{equation}
where $\mu,\nu\in\cH^{0,1}(X_{0},\End E_{0})$.

\begin{varthm} \label{2nd-variation}
Let $c_{1}(\lambda,\Vert\cdot\Vert_{Q})$ denote the first Chern form of the determinant line bundle $\lambda=\det\mathrm{ind}\,\delb\cong\det T^{\ast}\mathcal{N}$ with respect to Quillen's metric $\Vert\cdot\Vert_{Q}^{2}=\Vert\cdot\Vert^{2}(\det\Delta)^{-1}$. Then
$$c_{1}(\lambda,\Vert\cdot\Vert_{Q})=-\frac{1}{2\pi^{2}}\,\tilde{\Omega} +\delta,$$
where 
$$\delta=-\frac{2}{\pi}\sum_{i=1}^{n}\sum_{l,m=1}^{r_{i}}
\mathrm{sgn}(\alpha_{l}^{i}-\alpha_{m}^{i})(1 - 2|\alpha_{l}^{i}-\alpha_{m}^{i}|)k_{m}^{i}\Omega_{il}$$
is the cuspidal defect. 
\end{varthm}

\begin{proof} 
We need to prove that
\begin{gather*}
\delb_{\mathcal{N}}\del_{\mathcal{N}}\log\det\Delta=\Theta -\frac{\sqrt{-1}}{\pi}\,\tilde{\Omega} \\-4\sqrt{-1}\sum_{i=1}^{n}\sum_{l,m=1}^{r_{i}}
\mathrm{sgn}(\alpha_{l}^{i}-\alpha_{m}^{i})(1-2|\alpha_{l}^{i}-\alpha_{m}^{i}|)k_{m}^{i}\Omega_{il},
\end{gather*}
where the forms $\Theta$ and $\Omega_{il}$ were introduced in Sections 4.1 and 4.2 respectively.
It repeats almost verbatim the computation of 
$\displaystyle{\frac{\pa^{2}}{\pa\vep(\mu)\pa\overline{\vep(\nu)}}\log\det\Delta}$ 
in Theorem 2 of \cite{Takhtajan-Zograf89} using the variational formulas of Section 4.1. The only difference is in a
non-vanishing boundary term that appear after using \eqref{var-mu} and applying Stokes' theorem to the integral
$\int_{X_{0}}\ad L_{\bar\nu}\mu\wedge\mu$: 
\begin{align*}
\int\limits_{X_{0}}\ad L_{\bar\nu}\mu\wedge\psi & =-\int\limits_{X_{0}}\ad f_{\mu\bar{\nu}}\wedge\delb\psi +\lim_{Y\rightarrow\infty}\int\limits_{\pa F^{Y}}\ad f_{\mu\bar\nu}\wedge\psi.
\end{align*}
The first term in this formula is treated exactly as in \cite{Takhtajan-Zograf89}, whereas for the second term we get 
\begin{align*}
\lim_{Y\rightarrow\infty}\int\limits_{\pa F^{Y}}\ad f_{\mu\bar\nu}\wedge\psi
&=\lim_{Y\rightarrow\infty}\int\limits_{\pa F^{Y}}\tr(( f_{\mu\bar\nu}(z)\otimes I-I\otimes  f_{\mu\bar\nu}^{t}(z))\psi(z))dz \\
&=c_{1}+\dots +c_{n},
\end{align*}
where $c_{i}$ is the constant term of the Fourier expansion of 
$$\tr(( f_{\mu\bar\nu}(\sigma_{i}z)\otimes I-I\otimes  f_{\mu\bar\nu}^{t}(\sigma_{i}z))\psi(\sigma_{i}z))\sigma_{i}^{\prime}(z)$$ 
at the cusp $x_{i}$, $i=1,\dots,n$. From Lemmas \ref{psi-expansion} and \ref{f-mu-nu}, using the unitarity of $U_{i}$ and the definition of $\Omega_{il}$, we get
\begin{align*}
c_{i} & = \tr(( F_{\mu\bar\nu}^{i}\otimes I-I\otimes  (F^{i}_{\mu\bar\nu})^{t})C_{i}) \\
& = -4\sum_{l,m=1}^{r_{i}}
\mathrm{sgn}(\alpha_{l}^{i}-\alpha_{m}^{i})(1-2|\alpha_{l}^{i}-\alpha_{m}^{i}|)k_{m}^{i}\Omega_{il}.
\end{align*}
\end{proof}
\begin{remark}
By Lemma 4, the cuspidal defect can be rewritten as follows:
$$\delta = - \sum_{i=1}^{n}\sum_{l,m=1}^{r_{i}}\!
\mathrm{sgn}(\alpha_{l}^{i}-\alpha_{m}^{i})(1 - 2|\alpha_{l}^{i}-\alpha_{m}^{i}|)k_{m}^{i} c_{1}(\lambda_{il},\Vert\cdot\Vert_{il}).$$
\end{remark}

Now suppose that the parabolic structure is integral.
The elementary formula $\tr(\ad a\cdot\ad b)=2k\tr(ab)-2\tr a\tr b$, where $a,b\in\End\CC^{k}$, 
leads to the following result.

\begin{corollary} Let $c_{1}(\lambda_{0},\Vert\cdot\Vert_{Q})$ denote the first Chern form 
of the line bundle $\lambda_{0}\simeq T^{\ast}\mathcal{N}_{0}$ --- the restriction of the determinant line bundle $\lambda$ to the submanifold $\mathcal{N}_{0}\subset\mathcal{N}$. Then  
$$c_{1}(\lambda_{0},\Vert\cdot\Vert_{Q})=-\frac{k}{\pi^{2}}\,\Omega_{0}+\delta_{0},$$
where $\Omega_{0}$ and $\delta_{0}$ are the restrictions of the symplectic form $\Omega$ 
and the cuspidal defect $\delta$ to $\mathcal{N}_{0}$ respectively.
\end{corollary}

\subsection{Local index theorem in Atiyah-Singer's form}

Here we assume the existence of the universal endomorphism bundle $\cE\rightarrow X\times\mathcal{N}$, so that
$\cE|_{X\times\{E\}}\cong \End E$ for any point $\{E\}\in\mathcal{N}$. 
According to \cite{BY}, a universal bundle does exist for generic weight systems.
The metric $h_{\End E_0}$ defined fiberwise on the restriction
$\cE_0\rightarrow X_0\times\mathcal{N}$ (see Section 2.1) extends to a (pseudo)metric $h_{\cE}$ on $\cE$ that 
degenerates over the marked points $P_1,\dots,P_n\in X$. As in \cite{Takhtajan-Zograf89}, using the variational formulas from Section \ref{lie}, 
we can explicitly compute the curvature form related to $h_{\cE}$, and the corresponding Chern character form
$\mathrm{ch}(\cE,h_{\cE})$ is well defined as a current on $X\times\mathcal{N}$. As a result,
we can reformulate Theorem \ref{2nd-variation} in the Atiyah-Singer form with a cuspidal defect (cf. \cite{Takhtajan-Zograf91}). 

\begin{varthm}
Let $c_{1}(\lambda,\Vert\cdot\Vert_{Q})$ denote the first Chern form of the  determinant line bundle $\lambda=\det\mathrm{ind}\,\delb\cong\det T^{\ast}\mathcal{N}$ relative to Quillen's metric $\Vert\cdot\Vert_{Q}^{2}=\Vert\cdot\Vert^{2}(\det\Delta)^{-1}$. Then
$$c_{1}(\lambda,\Vert\cdot\Vert_{Q})=\pi_{\ast}(\mathrm{ch}_{2}(\cE,h_{\cE})) +\delta,$$
where $\mathrm{ch}_{2}(\cE,h_{\cE})$ is the $(2,2)$-component of the Chern character form of the universal 
endomorphism bundle $\cE$ relative to the metric $h_{\cE}$, $\pi_{\ast}:
C^{2,2}(X\times\mathcal{N})\rightarrow C^{1,1}(\mathcal{N})$ denotes integration along the fibers of the projection $\pi: X\times\mathcal{N}\rightarrow\mathcal{N}$, and 
$$\delta=-\frac{2}{\pi}\sum_{i=1}^{n}\sum_{l,m=1}^{r_{i}}
\mathrm{sgn}(\alpha_{l}^{i}-\alpha_{m}^{i})(1-2|\alpha_{l}^{i}-\alpha_{m}^{i}|)k_{m}^{i}\Omega_{il}$$
is the cuspidal defect. 
\end{varthm}

\begin{proof} Repeating the argument in \cite{Takhtajan-Zograf89}, 
it is not difficult to show that $\pi_{\ast}(\mathrm{ch}_{2}(\cE,h_{\cE}))=-\frac{1}{2\pi^{2}}\,\tilde{\Omega}$,
and the assertion immediately follows from Theorem 1.
\end{proof}

\subsection{A simple example}

Consider stable parabolic bundles of rank $2$ and parabolic degree $0$ with a single marked point. 
The parabolic structure is given by a complete flag $\CC^{2}\supset L\supset\{0\}$ at $P\in X$
(where $L$ is a line in $\CC^{2}$) with multiplicities
$0<\alpha_{1}<\alpha_{2}<1$. For an integral parabolic structure we have $\alpha_{1}+\alpha_{2}=1$, so that $\alpha_{1}=\alpha, 
\;\alpha_{2}=1-\alpha$, where $0<\alpha<\frac{1}{2}$. Such a parabolic structure is associated with an admissible 
$SU(2)$-representation $\rho$ of the fundamental group of $X\setminus P$, where 
the matrix $\rho(S)$ has eigenvalues $e^{2\pi\sqrt{-1}\alpha}$ and
$e^{-2\pi\sqrt{-1}\alpha}$. Without loss of generality we can assume that $\rho(S)$ is diagonal and the cusp lying over
the marked point $P$ is $\infty$. The cuspidal defect in this case is
$$\delta=-\frac{4(1-4\alpha)}{\pi}\Omega_{12} =-2(1-4\alpha)\,c_{1}(\lambda_{12}, \Vert\cdot\Vert_{12}),$$
where $\lambda_{12}$ is the line bundle on $\mathcal{N}$ with the Hermitian metric $\Vert\cdot\Vert_{12}$, defined in Section \ref{bundles}. 

In the simplest case of a pointed torus the moduli space $\mathcal{N}$ is just a complex projective
line, and $\lambda_{12}$ is the tautological line bundle on $\mathcal{N}\cong\mathbb{C}P^1$.
Since the bundles  $\lambda_0\cong T^{\ast}\mathbb{C}P^1$ and $\lambda_{12}$ have degrees $-2$ and $-1$ respectively, we get
\begin{align*}
\int\limits_{\mathcal{N}}c_{1}(\lambda_{0},\Vert\cdot\Vert_{Q})  =-2\quad\text{and}\quad
\int\limits_{\mathcal{N}}c_{1}(\lambda_{12},\Vert\cdot\Vert_{12})  =-1.
\end{align*}
By Corollary 1,
$$c_{1}(\lambda_{0},\Vert\cdot\Vert_{Q})= -\frac{2}{\pi^{2}}\,\Omega_{0} -2(1-4\alpha)c_{1}(\lambda_{12}, \Vert\cdot\Vert_{12}),$$
so that for the volume of $\mathcal{N}$ we have
\begin{equation} \label{volume}
Vol(\mathcal{N})=\int_{\mathcal{N}}\Omega_0 = 2\pi^{2}(1-2\alpha).
\end{equation}
Now let us recall a fascinating formula of Witten for the symplectic volume of the moduli space of parabolic $SU(2)$-bundles \cite[Formula 3.18]{Witten}. In our notation it reads
$$Vol(\mathcal{N})=2^{2g-1+n}\pi^{4g-4+n}\sum_{m=1}^\infty\frac{\prod_{i=1}^n\sin(2\pi\alpha_i m)}{m^{2g-2+n}},$$
where $g$ is the genus of the Riemann surface, $n$ is the number of marked points on it, and $\alpha_i$ are the weights at the marked points.
For the pointed torus it gives
$$Vol(\mathcal{N})=4\pi\sum_{m=1}^\infty\frac{\sin(2\pi\alpha m)}{m}=2\pi^{2}(1-2\alpha),$$
in agreement with our formula \eqref{volume}. In fact, Theorem 1 suggests an alternative method of computing symplectic volumes of moduli spaces of parabolic bundles (note that Witten' volume computation relies on the famous Verlinde formula).

\bibliographystyle{amsalpha}
\bibliography{vb}
\end{document}